\documentclass{article}      
\usepackage{graphicx}
\usepackage{xcolor}

\usepackage{amssymb}
\usepackage{eucal}
\usepackage{amsmath}
\usepackage{graphicx}
\usepackage{float}
\usepackage{framed}
\usepackage{hyperref}

\newcommand{\wt}{\widetilde}

\newcommand{\R}{{\mathbb  R}}  \numberwithin{equation}{section} \newtheorem{thm}{\bf Theorem}[section]
\newtheorem{lem}[thm]{\bf Lemma} \newtheorem{prop}[thm]{\bf Proposition} 
\newtheorem{cor}[thm]{\bf Corollary}
\newtheorem{defn}{\bf Definition}[section]

\DeclareMathOperator{\tr}{tr}

\begin{document}

\title{Optimization on the real symplectic group}

\author{Petre Birtea, Ioan Ca\c su, Dan Com\u{a}nescu\\
{\small Department of Mathematics, West University of Timi\c soara} 
\\
{\small Bd. V. P\^ arvan, No 4, 300223 Timi\c soara, Rom\^ania}\\
{\small Email: petre.birtea@e-uvt.ro, ioan.casu@e-uvt.ro, dan.comanescu@e-uvt.ro}}
\date{}

\maketitle

\begin{abstract}
We regard the real symplectic group $Sp(2n,\R)$ as a constraint submanifold of the $2n\times 2n$ real matrices $\mathcal{M}_{2n}(\R)$ endowed with the Euclidean (Frobenius) metric, respectively as a submanifold of the general linear group $Gl(2n,\R)$ endowed with the (left) invariant metric. For a cost function that defines an optimization problem on the real symplectic group we give a necessary and sufficient condition for critical points and we apply this condition to the particular case of a least square cost function. In order to characterize the critical points we give a formula for the Hessian of a cost function defined on the real symplectic group, with respect to both considered metrics. For a generalized Brockett cost function we present a necessary condition and a sufficient condition for local minimum. We construct a retraction map that allows us to detail the steepest descent and embedded Newton algorithms for solving an optimization problem on the real symplectic group.
\end{abstract}

\maketitle

\section{Introduction}

Recently there has been a great interest for optimization problems over smooth manifolds, which appear in the context of various applications, see \cite{absil}, \cite{fan}, \cite{fiori}. An important class of such problems are defined on the real symplectic group. For example, the minimization of the least square distance function on the symplectic group is extensively studied in \cite{wu-1} and \cite{wu-2}. This cost function has important applications in studying the fidelity  of dynamical gates in quantum analog computation and in the control of beam systems in particle accelerators, see \cite{wu-1}, \cite{mahony}, \cite{dragt}.

In Section \ref{First-order}, following  \cite{birtea-comanescu}, \cite{Birtea-Comanescu-Hessian}, \cite{5-electron}, we present a method to study the critical points and their nature for a cost function defined on a constraint manifold, which is embedded in a larger Riemannian manifold (usually an Euclidean manifold). 
The method consists in constructing the so-called \textbf{embedded gradient vector field} defined on the regular points of the ambient space with respect to the constraint functions. This vector field has the property that on the constraint manifold it coincides with the gradient of the cost function with respect to the induced metric. The Lagrange multiplier functions, which appear in the expression of the embedded gradient vector field, have the property to coincide with the classical Lagrange multipliers in critical points. We also emphasize an explicit formula for the Hessian of a cost function. Next, we describe the real symplectic group $Sp(2n,\R)$ as a constraint manifold.

In subsections \ref{euclidean-metric} and \ref{critical-euclidean} we first embed the symplectic group in the Euclidean space $\mathcal{M}_{2n}(\R)$, endowed with the Euclidean (Frobenius) metric. We compute the embedded gradient vector field with respect to this metric and we write it in matrix form. We also organize the Lagrange multiplier functions in matrix form and we show that this matrix verifies a Sylvester equation, which has a unique solution. We obtain a necessary and sufficient condition for critical points. We apply the result for a least square cost function and we recover the condition for critical points previously found in \cite{wu-2}.

In Subsection \ref{invariant-10} we embed the symplectic group in the general linear group $Gl(2n,\R)$, endowed with the (left) invariant metric. Using the connection between the gradients of a function with respect to the Frobenius and (left) invariant metrics we compute the embedded gradient vector field with respect to the invariant metric in a matrix form. In the case of the invariant metric, we obtain an explicit and simpler formula for the matrix of the Lagrange multiplier functions. 

Using the embedded gradient vector field with respect to the (left) invariant metric we detail the steepest descent algorithm on the real symplectic group. In order to do this, we also need to construct a retraction map. Using the idea of Cayley transform for the real symplectic group, see \cite{gosson} and \cite{machado-leite}, we introduce such a retraction map. We prove that this map verifies the conditions to be a retraction. All the steps of the algorithm are given in an explicit matrix form.

In Section \ref{second-order} we compute the Hessian matrix of a cost function defined on the real symplectic group in the previously considered metrics.
In Subsection \ref{euclidean-metric-second-order} we apply the general result in order to find a necessary condition and a sufficient condition for local minimum in the case of a generalized Brockett cost function. This cost function has been introduced in \cite{machado-leite} as an analogous of the classical Brockett cost function defined on the orthogonal group and introduced in \cite{brockett}.
In Subsection \ref{invariant-metric-second-order}, in order to obtain the Hessian matrix in the case of invariant metric, we first  compute the covariant derivative with respect to this metric and the Hessian of a smooth function defined on $Gl(2n,\R)$ applied on two left invariant vector fields. Using these formulas, we explicitly compute the Hessian matrix of the cost function with respect to the invariant metric.
Following the description of the embedded Newton algorithm on manifolds, as described in \cite{5-electron}, we detail its formulation on the real symplectic group. 

Geodesic-based numerical algorithms for optimization problems on the real symplectic group have been previously presented in \cite{fiori} and \cite{wang-sun-fiori}.

\section{First order optimality conditions on the real symplectic group}\label{First-order}

Let $\mathcal{S}\subset \mathfrak{M}$ be a submanifold of a Riemannian manifold $(\mathfrak{M},{\bf g})$, that can be described by a set of constraint functions, i.e. $\mathcal{S}={\bf F}^{-1}(c)$, where ${\bf F}=(F_1,\dots,F_k):\mathfrak{M}\rightarrow \R^k$ is a smooth map and $c\in \R^k$ is a regular value of ${\bf F}$. We endow $\mathcal{S}$ with the induced metric, hence $(\mathcal{S},{\bf g}_{_{ind}})$ becomes itself a Riemannian manifold. 

For solving optimization problems one needs, in general, to compute the gradient vector field and the Hessian operator of a smooth cost function $\widetilde{G}:(\mathcal{S},{\bf g}_{_{ind}})\rightarrow \R$. The Riemannian geometry of the submanifold $\mathcal{S}$ can be more complicated than the Riemannian geometry of the ambient manifold $\mathfrak{M}$. In what follows, we show how we can compute the gradient vector field and the Hessian operator of $\widetilde{G}$ using only the geometry of the ambient manifold $(\mathfrak{M},{\bf g})$.

Let $G:(\mathfrak{M},{\bf g})\rightarrow \R$ be a smooth prolongation of $\widetilde{G}$. In \cite{birtea-comanescu}, \cite{Birtea-Comanescu-Hessian}, \cite{5-electron}, it has been proved that
\begin{equation}\label{got}
\nabla _{{\bf g}_{_{ind}}}\widetilde{G}(s)=\partial_{\bf g} G(s),\,\,\,\forall s\in \mathcal{S},
\end{equation}
where $\partial_{\bf g} G$ is the unique vector field defined on the open set of regular points of the constraint function $\mathfrak{M}^{reg}\subset \mathfrak{M}$ that is tangent to the foliation generated by ${\bf F}$ having property \eqref{got}. We call $\partial_{\bf g} G$ {\bf the embedded gradient vector field} and it  is given by the following formula:
\begin{equation*}
\partial_{\bf g} G(s)=\nabla_{\bf g} G(s)-\sum\limits_{i=1}^k\sigma_{\bf g}^{i}(s)\nabla_{\bf g}F_i(s).
\end{equation*}
The Lagrange multiplier functions $\sigma_{\bf g}^{i}:\mathfrak{M}^{reg}\rightarrow \R$ are defined by the formula

\begin{equation}\label{sigma-101}
\sigma^i_{\bf g}(s):=\frac{\det \left(\text{Gram}_{(F_1,\ldots ,F_{i-1},G, F_{i+1},\dots,F_k)}^{(F_1,\ldots , F_{i-1},F_i, F_{i+1},...,F_k)}(s)\right)}{\det\left(\text{Gram}_{(F_1,\ldots ,F_k)}^{(F_1,\ldots ,F_k)}(s)\right)},
\end{equation}
where
\begin{equation*}\label{sigma}
\text{Gram}_{(g_1,...,g_s)}^{(f_1,...,f_r)}=\left[%
\begin{array}{cccc}
  {\bf g}(\nabla_{\bf g} g_1,\nabla_{\bf g}f_{1}) & ... & {\bf g}(\nabla_{\bf g} g_s,\nabla_{\bf g} f_{1}) \\
  \vdots & \ddots & \vdots \\
  {\bf g}(\nabla_{\bf g} g_1,\nabla_{\bf g}f_r) & ... & {\bf g}(\nabla_{\bf g} g_s,\nabla_{\bf g} f_r) 
\end{array}%
\right].
\end{equation*}

Also, in \cite{Birtea-Comanescu-Hessian}, \cite{5-electron} it has been proved that 
\begin{equation}\label{Hessian-general}
\text{Hess}_{{\bf g}_{_{ind}}}\, \widetilde{G}(s) =\left(\text{Hess}_{\bf g}\,G(s)-\sum_{i=1}^k\sigma_{\bf g}^{i}(s)\text{Hess}_{\bf g}\, F_i(s)\right)_{|T_s \mathcal{S}\times T_s \mathcal{S}}.
\end{equation}

In what follows the submanifold $\mathcal{S}$ is the symplectic group $Sp(2n,\R)$ and
let $\widetilde{G}:Sp(2n,\R)\rightarrow \R$ be a smooth cost function. We will endow the symplectic group with two different Riemannian metrics and solve the optimization problem generated by $\widetilde{G}$ in each case. 
More precisely, our optimization problem is
$$
\underset{M\in Sp(2n,\R)}{\mathrm{argmin}}~\widetilde{G}(M).
$$

The real symplectic group is defined by 
$$Sp(2n,\R)=\{M\in \mathcal{M}_{2n}(\R)\,|\,M^TJM=J\}, $$
where 
$$J=\left[
\begin{array}{cc}
\mathbb{O}_{n} & \mathbb{I}_n \\
-\mathbb{I}_n & \mathbb{O}_{n} 
\end{array}\right].
$$ 
A symplectic matrix can be written in the following block form:
$$M=\left[
\begin{array}{cc}
A & B \\
C & D 
\end{array}\right],$$
where $A,B,C,D\in \mathcal{M}_n(\R)$ verify the equalities
\begin{equation*}
A^TC=C^TA,\,\,\,B^TD=D^TB,\,\,\,A^TD-C^TB=\mathbb{I}_n.
\end{equation*}
We denote 
$$A=[ {\bf a}_1,...,{\bf a}_n],\,\,\,B=[ {\bf b}_1,...,{\bf b}_n],\,\,\,C=[ {\bf c}_1,...,{\bf c}_n],\,\,\,D=[ {\bf d}_1,...,{\bf d}_n],$$
with ${\bf a}_1,...,{\bf a}_n,{\bf b}_1,...,{\bf b}_n,{\bf c}_1,...,{\bf c}_n,{\bf d}_1,...,{\bf d}_n$ are column vectors in $\R^n$.

In this paper we embed the symplectic group in two larger Riemannian manifolds: $\mathcal{M}_{2n}(\R)$ endowed with the Euclidean (Frobenius) metric, respectively the general linear group $Gl(2n,\R)$ endowed with the invariant metric.
In the first case, the symplectic group is described by the following $2n^2-n$ functionally independent constraint functions 
${\bf F}^{AC}_{ij}, \,{\bf F}^{BD}_{ij},\,{\bf F}_{ij}:\mathcal{M}_{2n}(\R)\rightarrow \R$ given by
\begin{align}\label{constraints-1}
& {\bf F}^{AC}_{ij}(M)  = \left<{\bf a}_i,{\bf c}_j\right>-\left<{\bf a}_j,{\bf c}_i\right> ,\,\,\,1\leq i<j\leq n \nonumber\\
& {\bf F}^{BD}_{ij}(M)  = \left<{\bf b}_i,{\bf d}_j\right>-\left<{\bf b}_j,{\bf d}_i\right> ,\,\,\,1\leq i<j\leq n \\
& {\bf F}_{ij}(M)  = \left<{\bf a}_i,{\bf d}_j\right>-\left<{\bf c}_i,{\bf b}_j\right>-\delta_{ij} ,\,\,\,i,j\in\{1,2,...,n\}. \nonumber
\end{align}
In the second case, the constraint functions are the restrictions of $ {\bf F}^{AC}_{ij}$, $ {\bf F}^{BD}_{ij}$, and $ {\bf F}_{ij}$ to $Gl(2n,\R)$.


\subsection{The embedded gradient vector field using the Euclidean metric}\label{euclidean-metric}
 In the first case, we regard the symplectic group as an embedded manifold in the Euclidean space $(\mathcal{M}_{2n}(\R),\left<\cdot,\cdot\right>_{Euc})$, where $\left<X,Y\right>_{Euc}=\tr({X^TY})$ is the Frobenius scalar product. We prolong $\widetilde {G}$ to a smooth cost function $G:\mathcal{M}_{2n}(\R)\rightarrow \R$, i.e. $G_{|Sp(2n,\R)}=\widetilde{G}$. 


The embedded gradient vector field on $(\mathcal{M}_{2n}(\R),\left<\cdot,\cdot\right>_{Euc})$ associated to the constraint functions \eqref{constraints-1} is
\begin{equation}\label{partial-G-general}
\partial_{Euc}{G}=\nabla_{Euc} G-\sum_{1\leq i<j\leq n}\sigma^{ij}_{AC}\nabla_{Euc} {\bf F}^{AC}_{ij}-\sum_{1\leq i<j\leq n}\sigma^{ij}_{BD}\nabla_{Euc} {\bf F}^{BD}_{ij}-\sum_{i,j=1}^n\sigma^{ij}\nabla_{Euc} {\bf F}_{ij},
\end{equation}
where $\sigma_{AC}^{ij}, \sigma_{BD}^{ij}, \sigma^{ij}:(\mathcal{M}_{2n}(\R),\left<\cdot,\cdot\right>_{Euc})\rightarrow \R$ are the Lagrange multiplier functions \eqref{sigma-101} as introduced in \cite{birtea-comanescu}, \cite{Birtea-Comanescu-Hessian}.

Defining the matrices:
$$\Sigma_{AC}:=\left[
\begin{array}{cccc}
0 & \sigma_{AC}^{12} & ... & \sigma_{AC}^{1n}  \\
-\sigma_{AC}^{12} & 0 & ... & \sigma_{AC}^{2n}  \\
... & ... & ... & ... \\
-\sigma_{AC}^{1n} & -\sigma_{AC}^{2n} & ... & 0
\end{array}\right],
\Sigma_{BD}:=\left[
\begin{array}{cccc}
0 & \sigma_{BD}^{12} & ... & \sigma_{BD}^{1n}  \\
-\sigma_{BD}^{12} & 0 & ... & \sigma_{BD}^{2n}  \\
... & ... & ... & ... \\
-\sigma_{BD}^{1n} & -\sigma_{BD}^{2n} & ... & 0
\end{array}\right],
\Sigma:=\left[\sigma^{ij}\right]_{i,j=\overline{1,n}},
$$
by a straightforward computation we obtain the following matrix form for the embedded gradient vector field:
\begin{equation}\label{partial-Euc-matrix-form}
\partial_{Euc}{G}(M)=\nabla_{Euc} G(M)+JM\wt{\Sigma}(M),
\end{equation}
where we construct the $2n\times 2n$ skew-symmetric matrix $\wt{\Sigma}(M)$ as
$$\wt{\Sigma}(M)=\left[
\begin{array}{cc}
\Sigma_{AC}(M) & \Sigma(M) \\
-\Sigma(M)^T & \Sigma _{BD}(M)
\end{array}\right].$$

In what follows we denote by $\text{vec}(X)\in \R^{4n^2}$ the column vectorization of the matrix $X\in \mathcal{M}_{2n}(\R)$.  The matrix equality \eqref{partial-Euc-matrix-form} has been obtained by applying the operator $\text{vec}^{-1}$ to the vectorial equality \eqref{partial-G-general}.

Although the formulas \eqref{sigma-101}  are explicit, for the symplectic case they are extremely complicated from a computational point of view.  As an alternative approach, we search for an equation whose solution is $\wt{\Sigma}$.
By construction, see \cite{birtea-comanescu}, \cite{Birtea-Comanescu-Hessian}, $\partial_{Euc} G(M)\in T_{M}Sp(2n,\R)$ or equivalently, for every $M\in Sp(2n,\R)$ there exists a symmetric matrix $S(M)\in Sym(2n,\R)$, see \cite{Marsden-Ratiu}, such that  $\partial_{Euc} G(M)=MJS(M)$. Multiplying to the left with $(MJ)^{-1}=-M^TJ$, we obtain 
$$S(M)=-M^TJ\partial_{Euc} G(M).$$
The symmetry condition of $S(M)$ and the expression of $\partial_{Euc} G(M)$ lead to the following equation
$$-M^TJ(\nabla_{Euc} G(M)+JM\wt{\Sigma}(M))=(\nabla_{Euc} G(M)^T+\wt{\Sigma}(M)M^TJ)JM,$$
or equivalently 
\begin{equation}\label{ecuatie-sigma}
M^TM \wt{\Sigma}(M)+\wt{\Sigma}(M)M^TM=\nabla_{Euc} G(M)^TJM+M^TJ\nabla_{Euc} G(M).
\end{equation}
The above equation is a particular case of a Sylvester matrix equation, which has a unique solution since the matrices $M^TM$ and $-M^TM$ have no common eigenvalues as being positive definite matrix, respectively negative definite matrix ($M^TM$ is a Gramian matrix with the determinant equal 1). Explicitly, by using Kronecker notation, we have
\begin{equation*}
\text{vec}\,\wt{\Sigma}(M)=\left(\mathbb{I}_{2n}\otimes (M^TM)+(M^TM)\otimes \mathbb{I}_{2n}\right)^{-1}\text{vec}(\nabla_{Euc} G(M)^TJM+M^TJ\nabla_{Euc} G(M)).
\end{equation*}


By a straightforward computation we obtain the following result for the case $n=1$.

\begin{prop}
For $n=1$ the equation \eqref{ecuatie-sigma} has the unique solution
$$\widetilde{\Sigma}(M)=\frac{1}{\tr (M^TM)}\cdot \left(\nabla_{Euc} G(M)^TJM+M^TJ\nabla_{Euc} G(M)\right).$$
\end{prop}

\subsection{Critical points and applications using the Euclidean metric}\label{critical-euclidean}

The next result gives a necessary and sufficient condition for a symplectic matrix to be a critical point for the cost function $\wt{G}$.

\begin{thm}\label{critical}
A matrix $M\in Sp(2n,\R)$ is a critical point for the cost function $\wt{G}:Sp(2n,\R)\rightarrow \R$ if and only if 
$$JM^T\nabla_{Euc} G(M)=\nabla_{Euc} G(M)^TMJ.$$
Moreover, for critical points of the cost function $\wt{G}$, the matrix of Lagrange multipliers has the explicit form 
\begin{equation}\label{Sigma-critic}
\wt{\Sigma}(M)=JM^T\nabla_{Euc} G(M).
\end{equation}
\end{thm}

{\it Proof.}
A necessary and sufficient condition for critical points of the cost function $\wt{G}$ is $\partial_{Euc} G(M)={\bf 0}$, see \cite{birtea-comanescu}, \cite{Birtea-Comanescu-Hessian}, and \cite{5-electron}. The equality $\partial_{Euc} G(M)={\bf 0}$ is equivalent with 
$$\wt{\Sigma}(M)=JM^T\nabla_{Euc} G(M).$$

(i) First we prove the necessary condition. By definition $\wt{\Sigma}(M)$ is a skew-symmetric matrix, which implies the matrix equality stated in the theorem.

(ii) To prove that the equality of the theorem is also a sufficient condition, we check that the matrix $JM^T\nabla_{Euc} G(M)$ verifies the equation \eqref{ecuatie-sigma}. Indeed, we have
\begin{align*}
M^TM (JM^T\nabla_{Euc} G(M))+(JM^T\nabla_{Euc} G(M))M^TM & = M^TJ\nabla_{Euc} G(M)+ (JM^T\nabla_{Euc} G(M))M^TM \\
& = M^TJ\nabla_{Euc} G(M)+(\nabla_{Euc} G(M)^TMJ)M^TM \\
& = M^TJ\nabla_{Euc} G(M)+\nabla_{Euc} G(M)^TJM.
\end{align*}
From the uniqueness of the solution of equation \eqref{ecuatie-sigma} it follows that we have
$JM^T\nabla_{Euc} G(M)=\wt{\Sigma}(M)$, which implies that the matrix $M$ is a critical point for the cost function $\wt{G}$.
\rule{0.5em}{0.5em}

\subsubsection*{Least square cost function}

For the particular case of a least square cost function defined on the symplectic group,  we recover the necessary and sufficient condition from \cite{wu-1}, \cite{wu-2}.  
For a given symplectic matrix $W$ this cost function has the expression 
$$\wt{G}(M)=\|M-W\|_{Euc}^2=\text{tr}(M-W)^T(M-W).$$
By using the natural prolongation of $\wt{G}$, we obtain $\nabla_{Euc} G(M)=M-W$. Therefore, the necessary and sufficient condition of Theorem \ref{critical} becomes 
$$JM^T(M-W)=(M-W)^TMJ,$$
which is equivalent with the condition of  Theorem 3.1 in \cite{wu-2}: 
$$M^TM-(M^TM)^{-1}=M^TW-(M^TW)^{-1}.$$
In \cite{wu-2} it has been proved that the matrix  $W$ is a global minimum and all other critical points are saddles.

In \cite{fiori} a similar least square cost function is considered, without the assumption that the matrix $W$ is symplectic. This cost function originates from a mean problem of a data-set of symplectic matrices.

\subsection{The embedded gradient vector field using the invariant metric}\label{invariant-10}
 
In the second case, we regard the symplectic group as an embedded manifold in the Riemannian manifold  $(Gl(2n,\R),\left<\cdot,\cdot\right>_{Inv})$, where $\left<X_M,Y_M\right>_{Inv}=\left<M^{-1}X_M,M^{-1}Y_M\right>_{Euc}=\tr(X_M^TM^{-T}M^{-1}Y_M)$ is a left invariant  scalar product on the Lie group $Gl(2n,\R)$ and $X_M,Y_M\in T_M Gl(2n,\R)$. We prolong $\widetilde {G}$ to a smooth cost function $G:Gl(2n,\R)\rightarrow \R$, i.e. $G_{|Sp(2n,\R)}=\widetilde{G}$. 
  
  The embedded gradient vector field on $(Gl(2n,\R),\left<\cdot,\cdot\right>_{Inv})$ associated to the constraint functions \eqref{constraints-1} is
$$\partial_{Inv}{G}=\nabla_{Inv} G-\sum_{1\leq i<j\leq n}\gamma^{ij}_{AC}\nabla_{Inv} {\bf F}^{AC}_{ij}-\sum_{1\leq i<j\leq n}\gamma^{ij}_{BD}\nabla_{Inv} {\bf F}^{BD}_{ij}-\sum_{i,j=1}^n\gamma^{ij}\nabla_{Inv} {\bf F}_{ij},$$
where $\gamma_{AC}^{ij}, \gamma_{BD}^{ij}, \gamma^{ij}:(Gl(2n,\R),\left<\cdot,\cdot\right>_{Inv})\rightarrow \R$ are the Lagrange multiplier functions \eqref{sigma-101}.

We define the matrices:
$$\Gamma_{AC}:=\left[
\begin{array}{cccc}
0 & \gamma_{AC}^{12} & ... & \gamma_{AC}^{1n}  \\
-\gamma_{AC}^{12} & 0 & ... & \gamma_{AC}^{2n}  \\
... & ... & ... & ... \\
-\gamma_{AC}^{1n} & -\gamma_{AC}^{2n} & ... & 0
\end{array}\right],
\Gamma_{BD}:=\left[
\begin{array}{cccc}
0 & \gamma_{BD}^{12} & ... & \gamma_{BD}^{1n}  \\
-\gamma_{BD}^{12} & 0 & ... & \gamma_{BD}^{2n}  \\
... & ... & ... & ... \\
-\gamma_{BD}^{1n} & -\gamma_{BD}^{2n} & ... & 0
\end{array}\right],
\Gamma:=\left[\gamma^{ij}\right]_{i,j=\overline{1,n}}.
$$

To compute $\nabla_{Inv}G$,  we prove the following formula that relates $\nabla_{Inv}G$ with $\nabla_{Euc}G$.
\begin{lem}
Let $M\in Gl(2n,\R)$ and $f:Gl(2n,\R)\rightarrow \R$ a smooth function. Then we have 
\begin{equation*}\label{legatura-nabla-inv-euc}
\nabla_{Inv}f(M)=MM^T\nabla_{Euc}f(M).
\end{equation*}
\end{lem}

{\it Proof.}
For an arbitrary $Z_M\in T_MGl(2n,\R)$, we have 
$$\left<\nabla_{Inv}f(M),Z_M\right>_{Inv}=df(M)\cdot Z_M=\left<\nabla_{Euc}f(M),Z_M\right>_{Euc}, $$
which implies 
$$\tr \left(\nabla_{Inv}f(M)^TM^{-T}M^{-1}Z_M\right)=\tr \left(\nabla_{Euc}f(M)^TZ_M\right).$$
The tangent space $T_MGl(2n,\R)$ is given by all matrices of the form $MX$ with $X\in \mathcal{M}_{2n}(\R)$. For all $X\in \mathcal{M}_{2n}(\R)$ we have
$$\tr \left(\left(\nabla_{Inv}f(M)^TM^{-T}M^{-1}-\nabla_{Euc}f(M)^T\right)MX\right)=0,$$
which is equivalent with 
$$\left<M^T\left(\nabla_{Inv}f(M)^TM^{-T}M^{-1}-\nabla_{Euc}f(M)^T  \right)^T,X\right>_{Euc}=0.$$
The above equality proves the formula.
\rule{0.5em}{0.5em}
We construct the $2n\times 2n$ skew-symmetric matrix 
$$\wt{\Gamma}(M)=\left[
\begin{array}{cc}
\Gamma_{AC}(M) & \Gamma(M) \\
-\Gamma(M)^T & \Gamma_{BD}(M)
\end{array}\right].$$

The embedded gradient vector field in the case of the invariant metric becomes
{\small 
\begin{align*}
\partial_{Inv}G & =MM^T\nabla_{Euc} G-MM^T\left(\sum_{1\leq i<j\leq n}\gamma^{ij}_{AC}\nabla_{Euc} {\bf F}^{AC}_{ij}+\sum_{1\leq i<j\leq n}\gamma^{ij}_{BD}\nabla_{Euc} {\bf F}^{BD}_{ij}+\sum_{i,j=1}^n\gamma^{ij}\nabla_{Euc} {\bf F}_{ij}\right)\\
& =MM^T \nabla_{Euc} G +MM^T JM\widetilde{\Gamma}(M) = MM^T\nabla_{Euc} G+MJ\widetilde{\Gamma}(M).
\end{align*}}
Because $\partial_{Inv}G(M)\in T_M Sp(2n,\R)$ there exists the symmetric matrix $S(M)$ such that 
$$MM^T\nabla_{Euc} G(M)+MJ\widetilde{\Gamma}(M)=MJS(M).$$
Multiplying both size with $(MJ )^{-1}$ we obtain 
$$-JM^T\nabla _{Euc}G(M)+\widetilde{\Gamma}(M)=S(M).$$
Imposing the symmetry of $S(M)$ we obtain 
$$\widetilde{\Gamma}(M)=\frac{1}{2}\left(JM^T\nabla _{Euc}G(M)+\nabla _{Euc}G(M)^TMJ\right).$$

The following result has been obtained in \cite{fiori}, \cite{wang-sun-fiori} using another reasoning. 

\begin{thm} [\cite{fiori}, \cite{wang-sun-fiori}]
For $M\in Sp(2n,\R)$ the embedded gradient vector field in the case of the invariant metric has the formula
\begin{equation*}\label{partial-Inv}
\partial_{Inv}G(M)=\frac{1}{2}\left(MM^T\nabla_{Euc}G(M)+MJ\nabla_{Euc}G(M)^TMJ\right).
\end{equation*}
Moreover, the matrix of Lagrange multipliers $\widetilde{\Gamma}$ functions has the explicit formula 
\begin{equation}
\label{Lagrange-multipliers-Inv}
\widetilde{\Gamma}(M)=\frac{1}{2}\left(JM^T\nabla _{Euc}G(M)+\nabla _{Euc}G(M)^TMJ\right).\nonumber
\end{equation}
\end{thm}

Using the above result one can obtain a necessary and sufficient condition for critical points, which proves to be equivalent with the condition 
from Theorem \ref{critical}.

The advantage of using the vector field  $\partial_{Inv}G$ versus using the vector field $\partial_{Euc}G$ is that the former already is in explicit form, while for writing down the latter we need to solve in advance the Sylvester equation \eqref{ecuatie-sigma} that gives the matrix of Lagrange multipliers functions in the case of Euclidean metric.
This aspect can be useful for solving optimization problems on the symplectic group using first order numerical algorithms.

\subsubsection*{Steepest descent algorithm on the real symplectic group}

For implementing a steepest descent algorithm on a manifold one needs to map a vector from the tangent space into a point of the manifold. One way to do this is to use the notion of retraction, see \cite{absil-carte}. 

\begin{defn} A {\it retraction} on a manifold $\mathcal{S}$ is a smooth mapping $\mathcal{R}: T\mathcal{S} \rightarrow  \mathcal{S}$ with the following properties:
\begin{itemize}
\item[(i)] For all $x\in \mathcal{S}$ we have $\mathcal{R}_x({\bf 0}_x) = x$, where ${\bf 0}_x$ is the zero element of $T_x\mathcal{S}$ and $\mathcal{R}_x$ is the restriction of $\mathcal{R}$ to the tangent vector space $T_x\mathcal{S}$.
\item[(ii)] With the canonical identification $T_{{\bf 0}_x}T_x\mathcal{S}\simeq T_x\mathcal{S}$, $\mathcal{R}_x$ satisfies
$$D_{{\bf 0}_x}\mathcal{R}_x({\bf v}_x)={\bf v}_x,\,\,\,\forall~ {\bf v}_x\in T_x\mathcal{S}.$$
\end{itemize}
\end{defn}

\begin{prop}
The map  defined by 
\begin{equation}\label{retraction-expression}
\mathcal{R}(M,MJS)=-M(S+2J)^{-1}(S-2J),
\end{equation}
where $M\in Sp(2n,\R),\,S\in Sym(2n,\R),\,\emph{and}\,\det(S+2J)\neq 0$, 
is a retraction for the real symplectic group $Sp(2n,\R)$.

\end{prop}

{\it Proof.}
Following an idea from \cite{gosson} and \cite{machado-leite}, where a Cayley transform for the symplectic group is constructed, we are looking for a retraction of the form 
$$\mathcal{R}_M(MJS)=-M(S-\alpha J)^{-1}(S+\beta J),$$
where $\alpha,\beta\in \R$ and $S\in Sym(2n,\R)$.

To verify condition (i) of the above definition we have the computation
$$\mathcal{R}_M(\mathbb{O}_{2n})=-M(-\alpha J)^{-1}(\beta J)=-\frac{\beta}{\alpha}MJ^2=\frac{\beta}{\alpha}M,$$
from which we obtain the condition $\alpha=\beta$.

In order to verify condition (ii) we rewrite the map $\mathcal{R}_M$ in the equivalent form:
$$\mathcal{R}_M(MJS)=-M(-M^TJMJS-\alpha J)^{-1}(-M^TJMJS+\alpha J).$$
We make the notation $X:=MJS$ and consequently we obtain
$$\mathcal{R}_M(X) =-M(-M^TJX-\alpha J)^{-1}(-M^TJX+\alpha J).$$
Taking the differential of $\mathcal{R}_M$ with respect to the variable $X$ in the direction of a vector $MJ\wt{S}\in T_MSp(2n,\R)$ we have\footnote{We use the following formulas for matrix differentials: 
$D(U^{-1})=-U^{-1}(DU)U^{-1}$; $D(UV)=(DU)V+U(DV).$
}

\begin{align*}
D_X\mathcal{R}_M(MJ\wt{S})  = & M(-M^TJX-\alpha J)^{-1}(-M^TJ)MJ\wt{S}(-M^TJX-\alpha J)^{-1}(-M^TJX+\alpha J) \\
& -M(-M^TJX-\alpha J)^{-1}(-M^TJ)MJ\wt{S} \\
 = & M(-M^TJX-\alpha J)^{-1}\wt{S}(-M^TJX-\alpha J)^{-1}(-M^TJX+\alpha J) \\
 & -M(-M^TJX-\alpha J)^{-1}\wt{S}.
\end{align*}
It follows that 
\begin{align*}
D_{\mathbb{O}_{2n}}\mathcal{R}_M(MJ\wt{S})  = M(-\alpha J)^{-1}\wt{S}(-\alpha J)^{-1}(\alpha J)-M(-\alpha J)^{-1}\wt{S}=-\frac{2}{\alpha}MJ\wt{S}.
\end{align*}
The condition (ii) is verified if and only if $\alpha=-2$.

It remains to be checked that $\mathcal{R}_M(MJS)$ given by \eqref{retraction-expression} is a symplectic matrix, i.e. 
$$(\mathcal{R}_M(MJS))^TJ\mathcal{R}_M(MJS)=J.$$
We have the following equivalent equalities:
$$(S+2J)(S+2J)^{-T}M^TJM(S+2J)^{-1}(S-2J)=J\,\Leftrightarrow $$
$$(S+2J)(S+2J)^{-1}(-J)^{-1}(S+2J)^{-1}(S-2J)=J\,\Leftrightarrow $$
$$(S-2J)^{-1}(-J)^{-1}(S+2J)^{-1}=(S+2J)^{-1}(-J)^{-1}(S-2J)^{-1}\,\Leftrightarrow $$
$$(S+2J)(-J)(S-2J)=(S-2J)(-J)(S+2J),$$
which is obviously true. In the above computation we have used the fact that if $S+2J $ is invertible, then $S-2J$ is also invertible, since $S-2J=(S+2J)^T$.
\rule{0.5em}{0.5em}

The steepest descent algorithm for the real symplectic group has the following generic steps.

{\begin{framed}
\begin{itemize}

\item [1.] Consider a smooth prolongation $G:Gl(2n,\R)\rightarrow \R$ of the cost function  $\wt{G}:Sp(2n,\R)\rightarrow \R$.

\item [2.] Compute the matrix $\nabla_{Euc} G(M)$.

\item [3.] Input $M_0\in Sp(2n,\R)$ and $k=0$.

\item [4.] {\bf repeat}

$$M_{k+1}=-M_k(S_k+2J)^{-1}(S_k-2J).$$

{\bf until} $M_{k+1}$ sufficiently minimizes $\wt{G}$.

\medskip

{\bf Explanation}: the recurrence in the steepest descent algorithm on the real symplectic group is
$$M_{k+1}=\mathcal{R}_{M_k}(-\lambda_k\partial_{Inv}G(M_k)),$$
where $\lambda_k\in \R$ is a conveniently chosen length step (it has to verify $\det(S_k+2J)\neq 0$).

In order to use the explicit form of the retraction and taking into account the representation of the tangent vectors to the real symplectic group we need to solve the matrix equation $-\lambda_k\partial_{Inv} G(M_k)=M_kJS_k, $ with the unknown $S_k\in Sym(2n,\R)$.

This equation has the solution
$$S_k=\frac{\lambda_k}{2}\left(JM_k^T\nabla_{Euc} G(M_k)-(\nabla_{Euc} G(M_k))^TM_kJ\right).$$ 

\end{itemize}
\end{framed}

\section{Second order optimality on the real symplectic group}\label{second-order}

In this section we compute the Hessian of a cost function defined on the symplectic group with respect to the two metrics considered above. In the case of a generalized Brockett cost function, introduced in \cite{machado-leite},  we give a necessary condition and a sufficient condition for local minimum.

\subsection{Second order optimality using the Euclidean metric}\label{euclidean-metric-second-order}
In the first case, we regard the symplectic group as an embedded manifold in the Euclidean space $(\mathcal{M}_{2n}(\R),\left<\cdot,\cdot\right>_{Euc})$.
The formula \eqref{Hessian-general} for the Hessian matrix of the cost function $\widetilde{G}$, as given in \cite{Birtea-Comanescu-Hessian} and \cite{5-electron}, becomes in this case:
$$\text{Hess}_{Euc} \,\widetilde{G}(M):T_M Sp(2n,\R)\times T_M Sp(2n,\R)\rightarrow \R,$$
{\small \begin{align}\label{Hessian-formula-symplectic}
\text{Hess}_{Euc}\widetilde{G}(M) & =\left(\text{Hess}_{Euc}G(M)-\sum_{1\leq i<j\leq n}\sigma^{ij}_{AC}(M)\text{Hess}_{Euc} {\bf F}^{AC}_{ij}(M)\right. 
 \left.-\sum_{1\leq i<j\leq n}\sigma^{ij}_{BD}(M)\text{Hess}_{Euc}{\bf F}^{BD}_{ij}(M)\right. \nonumber \\
& \left. -\sum_{i,j=1}^n\sigma^{ij}(M)\text{Hess}_{Euc} {\bf F}_{ij}(M)\right)_{|T_M Sp(2n,\R)\times T_M Sp(2n,\R)}.
\end{align}}

By a straightforward computation we obtain the following formulas for the Hessian matrices of the constraint functions:
\begin{align}\label{Hess-euc-constraints}
 \text{Hess}_{Euc} {\bf F}^{AC}_{ij}(M) & =\Omega_{ij}^{AC}\otimes J,\,\,\,1\leq i<j\leq n, \nonumber \\
 \text{Hess}_{Euc}{\bf F}^{BD}_{ij}(M) & =\Omega_{ij}^{BD}\otimes J,\,\,\,1\leq i<j\leq n, \\
 \text{Hess}_{Euc}{\bf F}_{ij}(M) & =\Omega_{ij}\otimes J,\,\,\,1\leq i,j\leq n,\nonumber
\end{align}
where 
{\small $$\Omega_{ij}^{AC}=
\left[
\begin{array}{cc}
{\bf e}_i\otimes {\bf e}_j^T- {\bf e}_j\otimes {\bf e}_i^T & \mathbb{O}_n \\
\mathbb{O}_n & \mathbb{O}_n
\end{array}\right],\Omega_{ij}^{BD}=
\left[
\begin{array}{cc}
\mathbb{O}_n & \mathbb{O}_n \\
\mathbb{O}_n & {\bf e}_i\otimes {\bf e}_j^T- {\bf e}_j\otimes {\bf e}_i^T
\end{array}\right],
\Omega_{ij}=
\left[
\begin{array}{cc}
\mathbb{O}_n & {\bf e}_i\otimes {\bf e}_j^T \\
-{\bf e}^T_i\otimes {\bf e}_j & \mathbb{O}_n
\end{array}\right].\footnote{The vectors ${\bf e}_1$, ...   ,${\bf e}_n$ form the canonical basis in the Euclidean space $\R^n$.}
$$}
Substituting the above expressions in \eqref{Hessian-formula-symplectic} we obtain the following formula for the Hessian of a cost function defined on $Sp(2n,\R)$.

\begin{thm}\label{Formula-Hess-11} The Hessian of the cost function $\widetilde{G}:Sp(2n,\R)\rightarrow \R$ is given by
\begin{equation}\label{Hessian-symplectic-compact-102}
\emph{Hess}_{Euc}\widetilde{G}(M) =\left(\emph{Hess}_{Euc}G(M)-\widetilde{\Sigma}(M)\otimes J\right)_{|T_M Sp(2n,\R)\times T_M Sp(2n,\R)}.
\end{equation}
\end{thm}

As a consequence, we compute the formula of the Hessian of $\wt{G}$ applied on two tangent vectors from $ T_M Sp(2n,\R)$. 
\begin{cor}
Let $MJS_1,MJS_2\in T_M Sp(2n,\R)$, with $S_1,S_2\in Sym(2n,\R)$, be two tangent vectors. Then, for $M\in Sp(2n,\R)$
\begin{equation}\label{corolar-hessian}
\emph{Hess}_{Euc}\wt{G}(M)(MJS_1,MJS_2)=\emph{Hess}_{Euc}{G}(M)(MJS_1,MJS_2)+\tr \left(S_1JS_2\wt{\Sigma}(M)\right).
\end{equation}
In a critical point $M_c$ of $\wt{G}$ we have $\wt{\Sigma}(M_c)=JM_c^T\nabla_{Euc} G(M_c)$ and 
\begin{equation*}
 \emph{Hess}_{Euc}\wt{G}(M_c)(M_cJS_1,M_cJS_2)=\emph{Hess}_{Euc}{G}(M_c)(M_cJS_1,M_cJS_2)+\tr \left(S_1JS_2JM_c^T\nabla_{Euc} G(M_c)\right).
\end{equation*}
\end{cor}

{\it Proof.} By using the formula \eqref{Hessian-symplectic-compact-102} of Theorem \ref{Formula-Hess-11} we obtain\footnote{$\text{vec}(U)^T(Y\otimes X)\text{vec}(V)=\tr(U^TXVY^T)$.} 
\begin{align*}
-(\widetilde{\Sigma}(M)\otimes J)(MJS_1,MJS_2) & =-\text{vec}^T(MJS_1)(\widetilde{\Sigma}(M)\otimes J)\text{vec}(MJS_2) \\
& = -\tr\left((MJS_1)^TJ(MJS_2)\wt{\Sigma}^T(M)\right) 
=-\tr\left( S_1J(M^TJM)JS_2\wt{\Sigma}(M)\right) \\
& =\tr\left(S_1JS_2\wt{\Sigma}(M)\right).
\end{align*}
The second result from the above Corollary immediately follows from \eqref{Sigma-critic}.
\rule{0.5em}{0.5em}

\subsubsection*{Local extrema for a generalized Brockett cost function}

Following the ideas from \cite{brockett} for the orthogonal group acting on symmetric matrices, in \cite{machado-leite} an analogous theory is presented for $P$--orthogonal matrices\footnote{For an orthogonal matrix $P\in O(2n,\R)$ we say that a matrix $X\in \mathcal{M}_{2n}(\R)$ is $P$--orthogonal if $X^TPX=P$.}  acting on $P$--symmetric matrices. If $P=J$ we are in the case of the symplectic group, acting on $J$--symmetric matrices.

We recall the notions of $J$--transpose and $J$--symmetric matrices. For a matrix $L\in Gl(2n,\R)$ its $J$--transpose is by definition the matrix $L^J:=-JL^TJ$. The matrix $L$ is called $J$--symmetric if $L=L^J$ and we consider the real vector space of all $J$--symmetric matrices
$$\mathcal{I}=\{L\in Gl(2n,\R)~|~L=L^J\Leftrightarrow L^TJ=JL\}.$$
The symplectic group $Sp(2n,\R)$ acts naturally by conjugation on $\mathcal{I}$ and for a $J$--symmetric matrix $Q$ its orbit $\mathcal{O}_Q$ is the set 
$$\mathcal{O}_Q=\{M^{-1}QM|~M\in Sp(2n,\R)\}.$$
Further following the analogy with \cite{brockett} we have a natural pseudo--scalar product $\left<\cdot,\cdot\right>_J:\mathcal{I}\times \mathcal{I}\rightarrow \R$,
$$\left<X,Y\right>_J:=\tr(X^JY)=-\tr(JX^TJY).$$      
The optimization problem is to find the (pseudo)--distance from a given $J$--symmetric matrix $N$ to the orbit $\mathcal{O}_Q$ of a given $J$--symmetric matrix $Q$. It has been proved in \cite{machado-leite} that this problem is equivalent with the following optimization problem on the symplectic group:
$$\underset{M\in Sp(2n,\R)}{\mathrm{argmin}}~-\tr(M^{-1}QMN).$$    

We consider the so-called {\it generalized Brockett cost function} $\wt{G}:Sp(2n,\R)\rightarrow \R$ given by $$\wt{G}(M)=-\tr(M^{-1}QMN).$$
Since we have $M^{-1}=-JM^TJ$ we obtain that 
$$\wt{G}(M)=\tr(QMNJM^TJ).$$
We prolong $\wt{G}$ to the smooth function $G:\mathcal{M}_{2n}(\R)\rightarrow \R$, $G(M)=\tr(QMNJM^TJ).$ By a straightforward computation\footnote{
$\nabla_{Euc}\tr(XMYM^TZ)=X^TZ^TMY^T+ZXMY$ (see \cite{petersen-pedersen})} we get 
\begin{equation*}
\nabla_{Euc}G(M)= Q^TJMJN^T+JQMNJ=2JQMNJ.
\end{equation*}
Applying Theorem \ref{critical} we obtain that a matrix $M_c\in Sp(2n,\R)$ is a critical point for $\wt{G}$ if and only if $M_c^TJQM_cN$ is a skew-symmetric matrix, i.e.
$$M_c^TJQM_cN=N^TM_c^TQ^TJM_c.$$ 
This necessary and sufficient condition has been previously obtained in \cite{machado-leite}, in an equivalent form.

In order to compute the Hessian matrix in a critical point $M_c$ using formula \eqref{Hessian-symplectic-compact-102}, we first compute the matrix of Lagrange multipliers $\wt{\Sigma}(M_c)$ given by \eqref{Sigma-critic}:
$$\wt{\Sigma}(M_c)=JM^T_c\nabla_{Euc}G(M_c)=2JM_c^TJQM_cNJ.$$
In an arbitrary point $M\in \mathcal{M}_{2n}(\R)$ we have\footnote{
$\nabla_{Euc}(XMY)=Y^T\otimes X$ (see \cite{fackler})}
$$\text{Hess}_{Euc}G(M)=\nabla_{Euc}(2JQMNJ)=2(NJ)^T\otimes (JQ)=-2(JN^T)\otimes (JQ).$$
Consider now a critical point $M_c\in Sp(2n,\R)$ of $\wt{G}$ and an arbitrary tangent vector $M_cJS\in T_{M_c}Sp(2n,\R)$ determined by the symmetric matrix $S\in\mathcal{M}_{2n}(\R)$.\\
Next, we compute 
$\text{Hess}_{Euc}G(M_c)(M_cJS,M_cJS)$ and $(\wt{\Sigma}(M_c)\otimes J)(M_cJS,M_cJS)$ as follows:
\begin{align*}
\text{Hess}_{Euc}G(M_c)(M_cJS,M_cJS) & = -2((JN^T)\otimes (JQ))(M_cJS,M_cJS) =-2\tr(JM_c^TJQM_cJSNJS),
\end{align*}
respectively
\begin{align*}
(\wt{\Sigma}(M_c)\otimes J)(M_cJS,M_cJS) & = \tr((M_cJS)^TJM_cJS\wt{\Sigma}^T(M_c))=-\tr(SJS\wt{\Sigma}(M_c)) \\
& = -2\tr(SJSJM_c^TJQM_cNJ)  = -2\tr(JM_c^TJQM_cNJSJS).
\end{align*}
Using \eqref{Hessian-symplectic-compact-102} we obtain
\begin{equation*}
\text{Hess}_{Euc}\wt{G}(M_c)(M_cJS,M_cJS)=-2\tr(JM_c^TJQM_cJSNJS)+2\tr(JM_c^TJQM_cNJSJS).
\end{equation*}

We have the following necessary condition and sufficient condition for a local minimum of $\wt{G}$.

\begin{thm} Let $\wt{G}:Sp(2n,\R)\rightarrow \R$ given by $\wt{G}(M)=\tr(QMNJM^TJ)$.\\
(i) A symplectic matrix $M_c$ is a critical point for $\wt{G}$ if and only if $$M_c^TJQM_cN=N^TM_c^TQ^TJM_c.$$
(ii) If a critical point $M_c\in Sp(2n,\R)$ is a local minimum for $\wt{G}$, then the inequality 
$$\tr(JM_c^TJQM_cNJSJS)\geq \tr(JM_c^TJQM_cJSNJS)$$
holds for any symmetric matrix $S$.

A sufficient condition for a critical point $M_c\in Sp(2n,\R)$ to be local minimum is that the above condition is a  strict inequality for any non-zero symmetric matrix $S$. 
\end{thm}

\subsection{The Hessian matrix using the invariant metric}\label{invariant-metric-second-order}

On the Riemannian manifold  $(Gl(2n,\R),\left<\cdot,\cdot\right>_{Inv})$ the Hessian of a smooth function $f: Gl(2n,\R)\rightarrow \R$ is given by the following formula 
\begin{equation}\label{formula-Hessian-Riemannian}
\text{Hess}_{Inv}\, f(X,Y)=X(Yf)-(\nabla^{Inv}_X\, Y)f,
\end{equation}
where $X,Y\in \mathcal{X}(Gl(2n,\R))$ and $\nabla^{Inv}:\mathcal{X}(Gl(2n,\R))\times \mathcal{X}(Gl(2n,\R))\rightarrow \mathcal{X}(Gl(2n,\R))$ is the covariant derivative induced by the Riemannian metric $\left<\cdot,\cdot\right>_{Inv}$. 

Following an idea from \cite{dolcetti} we obtain the formula for the covariant derivative $\nabla^{Inv}$.

\begin{thm} On the Riemannian manifold  $(Gl(2n,\R),\left<\cdot,\cdot\right>_{Inv})$ we have following formulas:

\begin{itemize}
\item [(a)] for two vector 
fields $X,Y\in \mathcal{X}(Gl(2n,\R))$
\begin{align*}
(\nabla^{Inv}_{X}Y)_{_M} = & X_M(Y)-\frac{1}{2}\cdot\left(Y_MM^{-1}X_M+X_MM^{-1}Y_M+MY_M^TM^{-T}M^{-1}X_M \right. \\
& \left. + MX_M^TM^{-T}M^{-1}Y_M -Y_MX_M^TM^{-T}-X_MY_M^TM^{-T}\right);
\end{align*}

\item [(b)] for two left invariant vector 
fields $X,Y\in \mathcal{X}(Gl(2n,\R))$, i.e. $X_M=MX_0$, $Y_M=MY_0$ with $X_0,Y_0\in \mathcal{M}_{2n}(\R)$, we have\footnote{We denote by $[X,Y]=XY-YX$ the bracket of the matrices $X,Y$. }
\begin{equation}\label{derivative-covariant-Inv}
(\nabla^{Inv}_{X}Y)_{_M}=\frac{1}{2}M\left([X_0,Y_0]-[Y_0^T,X_0]-[X_0^T,Y_0]\right).
\end{equation}
\end{itemize}
\end{thm}

The corresponding result for the the right invariant metric on $Gl(2n,\R)$  has been proved before in \cite{absil} using Koszul's formula. 

\begin{thm}
For a smooth function $f$ defined on the Riemannian manifold  $(Gl(2n,\R),\left<\cdot,\cdot\right>_{Inv})$ and for  two left invariant vector 
fields $X,Y\in \mathcal{X}(Gl(2n,\R))$, i.e. $X_M=MX_0$, $Y_M=MY_0$ with $X_0,Y_0\in \mathcal{M}_{2n}(\R)$, we have the following formula:
\begin{equation}\label{Hess-inv-Gl}
\emph{Hess}_{Inv}\, f(MX_0,MY_0)=\emph{Hess}_{Euc}\, f(MX_0,MY_0)+\left<\nabla_{Euc}f(M),M{Z}_{(X_0,Y_0)}\right>_{Euc},
\end{equation}
where 
\begin{equation*}\label{Z-0}
Z_{(X_0,Y_0)}=\frac{1}{2}\left(X_0Y_0+Y_0X_0+[Y_0^T,X_0]+[X_0^T,Y_0]\right).
\end{equation*}
\end{thm}

{\it Proof.}
Substituting \eqref{derivative-covariant-Inv} into 
\eqref{formula-Hessian-Riemannian} we have the following straightforward computations.
\begin{align*}
\text{Hess}_{Inv}\, f(MX_0,MY_0)  & = MX_0\left(D_Mf\cdot MY_0\right)-(\nabla^{Inv}_{MX_0}\, MY_0)f \\
& = \text{Hess}_{Euc}\, f(MX_0,MY_0)+D_Mf\cdot MX_0Y_0 \\
& -\frac{1}{2}D_Mf\left(MX_0Y_0-MY_0X_0-[Y_0^T,X_0]-[X_0^T,Y_0]\right) \\
&  =\text{Hess}_{Euc}\, f(MX_0,MY_0) \\
& +\left<\nabla_{Euc}f(M),\frac{1}{2}M\left(X_0Y_0+Y_0X_0+[Y_0^T,X_0]+[X_0^T,Y_0]\right)\right>_{Euc}.
\end{align*}
\rule{0.5em}{0.5em}

The following  result provides the formula for the Hessian of a cost function defined on the symplectic group with respect to the invariant metric.

\begin{thm}\label{hess-inv-th}
Let  $X_M,Y_M\in T_M Sp(2n,\R)$ with $X_M=MX_0=MJS_X$ and $Y_M=MY_0=MJS_Y$, where $S_X,S_Y\in Sym(2n,\R)$.
The Hessian of a cost function $\widetilde{G}:Sp(2n,\R)\rightarrow \R$ is given by
\begin{align*}
\emph{Hess}_{Inv}\,\wt{G}(X_M,Y_M)&=\left(\emph{Hess}_{Euc}\,{G}-\wt{\Gamma}(M)\otimes J\right)(X_M,Y_M)+\nonumber\\
&+\left<\nabla_{Euc} G(M)+JM\wt{\Gamma}(M),MZ_{(X_0,Y_0)}\right>_{Euc},
\end{align*}
where 
$$\widetilde{\Gamma}(M)=\frac{1}{2}\left(JM^T\nabla _{Euc}G(M)+\nabla _{Euc}G(M)^TMJ\right)$$
and 
$$Z_{(X_0,Y_0)}=\frac{1}{2}\left(X_0Y_0+Y_0X_0+[Y_0^T,X_0]+[X_0^T,Y_0]\right).$$
\end{thm}

{\it Proof.}
Using  \eqref{Hessian-general} and \eqref{Hess-inv-Gl} for the case of the symplectic group we have 
{ \begin{align*}\label{Hessian-formula-symplectic-invariant}
\text{Hess}_{Inv}\widetilde{G}(X_M,Y_M) & =\text{Hess}_{Inv}G(X_M,Y_M)-\sum_{1\leq i<j\leq n}\gamma^{ij}_{AC}(M)\text{Hess}_{Inv} {\bf F}^{AC}_{ij}(X_M,Y_M) \\
 & -\sum_{1\leq i<j\leq n}\gamma^{ij}_{BD}(M)\text{Hess}_{Inv}{\bf F}^{BD}_{ij}(X_M,Y_M)
  -\sum_{i,j=1}^n\gamma^{ij}(M)\text{Hess}_{Inv} {\bf F}_{ij}(X_M,Y_M) \\
  & =  \text{Hess}_{Euc}\, G(X_M,Y_M)+\left<\nabla_{Euc}G(M),M{Z}_{(X_0,Y_0)}\right>_{Euc} \\
  & - \sum_{1\leq i<j\leq n}\gamma^{ij}_{AC}(M)\left(   \text{Hess}_{Euc}\, {\bf F}^{AC}_{ij}(X_M,Y_M)+\left<\nabla_{Euc}{\bf F}^{AC}_{ij}(M),M{Z}_{(X_0,Y_0)}\right>_{Euc}  \right) \\
  & - \sum_{1\leq i<j\leq n}\gamma^{ij}_{BD}(M)\left(   \text{Hess}_{Euc}\, {\bf F}^{BD}_{ij}(X_M,Y_M)+\left<\nabla_{Euc}{\bf F}^{BD}_{ij}(M),M{Z}_{(X_0,Y_0)}\right>_{Euc}  \right) \\
  & - \sum_{i,j=1}^n\gamma^{ij}(M)\left(   \text{Hess}_{Euc}\, {\bf F}_{ij}(X_M,Y_M)+\left<\nabla_{Euc}{\bf F}_{ij}(M),M{Z}_{(X_0,Y_0)}\right>_{Euc}  \right).
\end{align*}}
Using the formulas \eqref{Hess-euc-constraints} and analogous arguments as in Theorem \ref{Formula-Hess-11} and equation \eqref{partial-Euc-matrix-form} we further obtain
\begin{align*}
\text{Hess}_{Inv}\widetilde{G}(X_M,Y_M) & = \left(\text{Hess}_{Euc}G-\widetilde{\Gamma}(M)\otimes J\right)(X_M,Y_M) 
 + \left<\nabla_{Euc} G(M)+JM\wt{\Gamma}(M),MZ_{(X_0,Y_0)}\right>_{Euc}.
\end{align*}
\rule{0.5em}{0.5em}

We notice that the Hessian matrix $\text{Hess}_{Inv}\widetilde{G}$ is given by an explicit formula, while the expression of the Hessian matrix $\text{Hess}_{Euc}\widetilde{G}$ contains the Lagrange multipliers matrix $\wt{\Sigma}$ that is given by the Sylvester matrix equation \eqref{ecuatie-sigma}. This aspect can be advantageous for solving optimization problems on the symplectic group using second order numerical algorithms.

\subsubsection*{Newton algorithm on the real symplectic group}

In what follows we adapt the generic embedded Newton algorithm on manifolds, as presented in \cite{5-electron}, to the specific case of the real symplectic group.

\begin{framed}
	\begin{itemize}
		\item [1.] Consider a smooth prolongation $G:Gl(2n,\R)\rightarrow \R$ of the cost function  $\wt{G}:Sp(2n,\R)\rightarrow \R$.
		
		\item [2.] Construct the constant matrices ${\bf e}_{(i,j)}=\frac{1}{2}J({\bf f}_i\otimes {\bf f}_j+{\bf f}_j\otimes {\bf f}_i)$, $1\leq i\leq j\leq 2n$, where the vectors ${\bf f}_1$, \dots   ,${\bf f}_{2n}$ denote the canonical basis in the Euclidean space $\R^{2n}$. The matrices $M{\bf e}_{(i,j)}$, $1\leq i\leq j\leq 2n$, form a basis for the tangent vector space $T_MSp(2n,\R)$.
		
		\item [3.] Compute the coordinate functions 
		\begin{equation*}\label{coefficients-g}
		g_{(i,j)}(M)=dG(M)\cdot M{\bf e}_{(i,j)},\,1\leq i\leq j\leq 2n.
		\end{equation*}
		
		\item [4.] Compute the matrices $\text{Hess}_{Euc}\,{G}(M)$, $\nabla _{Euc}G(M)$, $\widetilde{\Gamma}(M)$, and $Z_{\left({\bf e}_{(i,j)},{\bf e}_{(p,q)}\right)},$\\ for all $1\leq i\leq j\leq 2n,\,1\leq p\leq q\leq 2n$.
		
		\item [5.] Compute the components of the Hessian matrix $\text{Hess}_{Inv}\,\wt{G}$ of the cost function $\wt{G}$
		\begin{equation*}\label{HessR}
		h_{(i,j),(p,q)}(M)= \text{Hess}_{Inv}\,\wt{G}\left(M{\bf e}_{(i,j)},M{\bf e}_{(p,q)}\right),\,1\leq i\leq j\leq 2n,\,1\leq p\leq q\leq 2n,
		\end{equation*}
		using the formulas from Theorem \ref{hess-inv-th}.
		
		\item [6.] Input $M_0\in Sp(2n,\R)$ and $k=0$.
		
		\item [7.] {\bf repeat}
		
		$\bullet$ Solve the linear system with the unknowns  $v^{(i,j)}_{M_k}$, $1\leq i\leq j\leq 2n$
		\begin{equation*}\label{newton-equation}
		\sum_{(i,j)} h_{(i,j),(p,q)}(M_k)v^{(i,j)}_{M_k}=-g_{(p,q)}(M_k),\,1\leq p\leq q\leq 2n,
		\end{equation*}
which represents the Newton equation.
		
		$\bullet$ Construct the line search tangent vector $$ {\bf v}_{M_k}=\sum_{(i,j)} v^{(i,j)}_{M_k}M_k{\bf e}_{(i,j)}.$$
		
		$\bullet$ Compute the symmetric matrix $S_k=-\lambda_kM_k^TJ{\bf v}_{M_k}$, where $\lambda_k\in \R$ is a conveniently chosen length step. 
		
		$\bullet$ Compute $$M_{k+1}=-M_k(S_k+2J)^{-1}(S_k-2J)$$

{\bf until} $M_{k+1}$ sufficiently minimizes $\wt{G}$.
\medskip

{\bf Explanation}: the recurrence in the Newton algorithm on the real symplectic group is
$$M_{k+1}=\mathcal{R}_{M_k}(\lambda_k{\bf v}_{M_k}).$$

In order to use the explicit form of the retraction \eqref{retraction-expression} and taking into account the representation of the tangent vectors to the real symplectic group we need to solve the matrix equation $\lambda_k{\bf v}_{M_k}=M_kJS_k $ with the unknown $S_k\in Sym(2n,\R)$.
This equation has the solution $S_k=-\lambda_kM_k^TJ{\bf v}_{M_k}$.
		
	\end{itemize}
\end{framed}

\subsection*{Acknowledgements}
This work was supported by a grant of Ministery of Research and Innovation, CNCS-UEFISCDI, project number PN-III-P4-ID-PCE-2016-0165, within PNCDI III.

\end{document}